\documentclass[a4paper,11pt]{amsart}
\usepackage{amssymb}
\usepackage[english]{babel}
\usepackage[latin1]{inputenc}
\usepackage{graphicx}

\theoremstyle{plain}
\newcommand{\ep}{\varepsilon}
\newcommand{\R}{\mathbb R}

\begin{document}
\title{Convexity theory for the term structure equation}
\author{Erik Ekström and Johan Tysk}
\email{ekstrom@math.uu.se, johan.tysk@math.uu.se}
\subjclass[2000]{Primary 91B28; Secondary 35K99} \keywords{Interest
rate theory, bond options, convexity, parameter monotonicity,
log-convexity, log-concavity, affine term structure.}
\address{Department of
Mathematics, Uppsala University, Box 480, 75106 Uppsala, Sweden}

\newtheorem{theorem}{Theorem}[section]
\newtheorem{lemma}[theorem]{Lemma}
\newtheorem{corollary}[theorem]{Corollary}
\newtheorem{proposition}[theorem]{Proposition}
\newtheorem{definition}[theorem]{Definition}
\newtheorem{hypothesis}[theorem]{Hypothesis}
\newtheorem{assumption}[theorem]{Assumption}

\newenvironment{example}[1][Example]{\begin{trivlist}
\item[\hskip \labelsep {\bf Example}]}{\end{trivlist}}
\newenvironment{remark}[1][Remark]{\begin{trivlist}
\item[\hskip \labelsep {\bf Remark}]}{\end{trivlist}}

\begin{abstract}
We study convexity and monotonicity properties for prices 
of bonds and bond options when the short rate is modeled 
by a diffusion process. We provide conditions under which
convexity of the price in the short rate is guaranteed. 
Under these conditions the price is decreasing in the
drift and increasing in the volatility of the short rate.
We also study convexity properties of the logarithm of 
the price.
\end{abstract}

\maketitle

\section{Introduction}

Already in the seminal paper \cite{M}, convexity of the option price
in the underlying asset is discussed. In the last decade this issue
has attracted renewed interest in the literature, compare
\cite{A1}, \cite{BJ}, \cite{BR}, \cite{BGW}, \cite{E}, \cite{EJT}, \cite{ET1},
\cite{ET2}, \cite{EKJPS}, \cite{ER}, \cite{H} and \cite{JT1}. Given
a convex pay-off function one asks for what models this convexity is
preserved in the sense that the price also is a convex function of
the underlying asset at any fixed time prior to maturity. This
question is studied for diffusion models, for models with jumps and
for various option types including options written on several
underlying assets. The interest in convexity has at least three
reasons. Firstly, convexity is a fundamental qualitative property of
option prices (an even more fundamental qualitative property would
be monotonicity in the underlying asset provided the pay-off is
monotone, but such properties can usually be derived immediately).
Secondly, if the price is convex then it is also typically
increasing in the volatility, and in the case of jump-diffusion
models, also in the jump parameters. Thirdly, if a delta-hedger uses
a model that overestimates the true volatility, he or she will
obtain  a superhedge for the claim provided the price is convex.

Our aim with the present paper is to continue this study to bonds
and bond options, for which we regard the short rate as the
underlying process. Thus we study preservation of convexity for the
term structure equation instead of variants of the Black-Scholes
equation as is the case in the references above. Surprisingly little
has been done in this direction, with \cite{A} as a notable
exception. A motivation for this is perhaps that the third reason
for studying convexity mentioned above is not directly applicable
since the short rate is not a traded asset. However, the first two
stated motives remain valid also in interest rate theory.

We assume that the short rate is modeled under some given risk neutral probability
measure as a stochastic process $X=(X_t)_{t\geq 0}$ with dynamics
\begin{equation*}
dX_t=\beta(X_t,t)\,dt+\sigma(X_t,t)\,dB_t,
\end{equation*}
where $B$ is a standard Brownian motion and $\beta$ and $\sigma$ are
given functions of time and the current short rate. We first investigate
the convexity properties of the $T$-bond price
\begin{equation}
\label{u}
u(x,t)=E_{x,t}\left[\exp\left\{-\int_t^TX_s\,ds\right\}\right],
\end{equation}
where the indices indicate that $X_t=x$. Using the Feynman-Kac
theorem it follows that the bond price $u$ satisfies the term
structure equation
\begin{equation*}\label{termstructure}\left\{\begin{array}{ll}
u_t+\frac{\sigma^2}{2} u_{xx}+\beta u_x-xu=0 & \mbox{for }t<T\\
u=1& \mbox{for }t=T.\end{array}\right. \end{equation*} 
This enables us to use a mixture of stochastic techniques and methods from the theory of parabolic partial differential equations.

One should note that our pay-off function in the case of bonds is
identically equal to one, so it is both convex and concave. However,
convexity in $x$ is the natural property to consider: since the bond
price declines with $x$, convexity means that the absolute value of
this decline decreases with $x$. We certainly expect bond prices to
suffer a smaller decline if short rates move from $5\%$ to $6\%$
than if the short rates move from $1\%$ to $2\%$! In Section
\ref{convexity} we find conditions on the model for $X$ that
guarantee that convexity is preserved. More precisely, if
$\beta_{xx}\leq 2$ (in the sense of distributions if $\beta$ is not
twice continuously differentiable), then the model is convexity preserving. To our knowledge,
this condition is indeed fulfilled for all models of the short rate
that are used in practice, compare Table~\ref{table1} below. Moreover, the condition is sharp in the
sense that it is also a necessary condition for preservation of
convexity provided the coefficients of the model are regular enough,
see Theorem~\ref{nec}. Using the pathwise non-crossing property of
diffusions, it is easily seen that the bond price $u$
is decreasing in the drift $\beta$. In Section~\ref{mon} we 
show that the bond price in a convexity preserving model
is also increasing in the volatility $\sigma$. Thus the general 
relationship between convexity and monotonicity in the 
volatility known from option pricing theory extends to
interest rate theory.

Our study of convexity and monotonicity properties is formulated for
prices of options written with the short rate as the underlying
asset. Thus we study the function
\begin{equation}
\label{U}
U(x,t)=E_{x,t}\left[\exp\left\{-\int_t^TX_s\,ds\right\} g(X_T)\right],
\end{equation}
where $g$ is a given convex pay-off function (note that the case
$g\equiv 1$ corresponds to a bond). The main reason for extending
the study from bonds to options on the short rate is to be able to
study bond option prices, i.e. prices of options written on 
a bond price as the underlying asset. In fact, our general 
convexity and
monotonicity results allow us in Section~\ref{bond options} to
deduce properties of certain bond options. For example, 
the price of a
bond call option is convex in $x$ and therefore also increasing
in $\sigma$ for convexity preserving models.

It is also natural to consider convexity properties of the {\em logarithm} of the bond price. This is connected to the notion
of {\em duration}, i.e. the negative of the derivative of
the logarithm. The analogous concept for stock options is 
elasticity, compare \cite{B1}, \cite{B2},
\cite{ET3} and \cite{K}. We say that the price is log-convex if the
logarithm of the price is convex in $x$ and analogously for
log-concavity. Again, since the pay-off is constant for the bond it
is both log-convex and log-concave. Unlike the case of convexity,
however, both of these cases deserve consideration. 
According to the discussion above, convexity of a bond price means that the {\em absolute} value of the decline is decreasing in $x$. In contrast to
this, log-convexity means that the {\em relative} decline diminishes
when $x$ grows (a declining duration), and log-concavity 
means that the relative decline of the price increases 
(an increasing duration). In Section~\ref{Log-convexity} we 
show that if the drift $\beta$ is spatially concave and 
the diffusion coefficient $\sigma^2$ is spatially convex,
then log-convexity is preserved. Similarly, in
Section~\ref{Log-concavity} we show that 
log-concavity is preserved provided $\beta$ is convex and 
$\sigma^2$ is concave. If we insist that the model should
preserve both log-convexity and log-concavity, we arrive at 
models where the logarithm of the bond price is both convex 
and concave, i.e. linear. These are, of course, the models 
that admit an affine term structure (in our context these 
bond prices would be referred to as being log-affine).
Apart from admitting explicit bond prices, affine models 
play an important r\^ole in interest rate theory, compare 
for example \cite{BS}. Thus we recover the well-known sufficient condition that
$\beta$ and $\sigma^2$ are affine for the existence of an
affine term structure.

In the next section we present the assumptions on the
model parameters under which our results are presented.
If additional regularity of the coefficients is assumed,
bounds on the spatial derivatives of bond and option prices
can be obtained, see Section~\ref{aux}. In 
Section~\ref{continuity} a continuity result is provided,
thus allowing us to assume that the coefficients are 
regular, and the general results follow by approximating
the coefficients. As mentioned above, our main results are 
presented in Sections~\ref{convexity}-\ref{Log-concavity}, 
and a summary of these results is presented in 
Section~\ref{summary}.

One important feature of our approach is 
that it works just as well for models where $X$ 
never reaches zero, for models where the rate
can reach zero, as for models that allow negative 
interest rates. We thus avoid a case by case study. 
However, we will throughout the paper point out to 
which of the commonly studied models presented in
Table~\ref{table1} our results are applicable. For a more detailed
discussion of these models, see for example \cite{BM}.

\begin{table}[ht!]
\begin{center}
\begin{tabular}{|c|c|c|c|c|}
\hline Model & Dynamics & $X>0$ & AB & AO \\
\hline \hline V & $dX=k(\theta-X)\,dt+\sigma\,dB$ & No & Yes & Yes\\
\hline CIR & $dX=k(\theta-X)\,dt+\sigma\sqrt X\,dB$ & Yes &Yes & Yes\\
\hline D & $dX=bX\,dt+\sigma X\,dB$ & Yes & Yes & No\\
\hline EV & $dX=X(\eta-a\ln X)\,dt+\sigma X\,dB$ & Yes & No & No\\
\hline HW & $dX=k(\theta_t-X)\,dt+\sigma\,dB$ & No & Yes & Yes\\
\hline BK & $dX=X(\eta_t-a\ln X)\,dt+\sigma X\,dB$ & Yes & No & No\\
\hline MM & $dX=X(\eta_t-(\lambda-\frac{\gamma}{1+\gamma t})\ln
X)\,dt+\sigma X\,dB$ & Yes & No & No\\
\hline
\end{tabular}\caption{\label{table1} Some examples of short
rate models. The table is copied from \cite{BM}. AB stands for
analytic bond prices, and AO stands for analytic bond option prices.
When referring to these models below we will assume that $k$,
$\theta$, $a$, $\eta$ and $\gamma$ are positive and that
$\lambda\geq \gamma$.}
\end{center}
\end{table}

\section{Assumptions}
As explained in the introduction, we assume that the short rate $X$ is modeled as
\begin{equation}
\label{X} dX_t=\beta(X_t,t)\,dt+\sigma(X_t,t)\,dB_t
\end{equation}
for some Brownian motion $B$. Note that some of the short rate
models in Table~\ref{table1} are specified so that the interest rate
cannot fall below zero, whereas some models allow for negative
interest rates with positive probability. To unify the analysis, we
view the short rate process $X$ in \eqref{X} as specified on the
whole real line $\R$ with $\sigma$ and $\beta$ suitably extended
(for example to be 0) for negative values in case $X$ is specified
initially only on the positive real axis.

Throughout the article we make the following regularity and growth
assumptions.

\begin{assumption}
\label{holder} The functions $\beta,\sigma:
\R\times[0,\infty)\to\R$
are continuous functions, $\beta$ is locally Lipschitz 
continuous in
the $x$-variable, and $\sigma$ is locally H\"older(1/2) in the
$x$-variable. Moreover, there exists a constant $D$ such that
\begin{equation}
\label{growthsigma} \vert\sigma(x,t)\vert \leq
D(1+x^+)\end{equation} and \begin{equation} \label{growthbeta}
\vert\beta(x,t)\vert\leq D(1+\vert x\vert)\end{equation} for all
$(x,t)\in\R\times [0,\infty)$.
\end{assumption}

The conditions on $\sigma$ and $\beta$ guarantee a non-exploding
unique strong solution to \eqref{X} for any initial point $(x,t)\in
\R\times[0,\infty)$. The condition that $\sigma$ is bounded for
negative $x$ implies that, for any $T$, the price $u$ in \eqref{u}
of a $T$-bond is finite at all points $(x,t)\in\R\times [0,T]$,
see Corollary~\ref{cor} below. Without this condition, the bond
price $u$ may be infinite. Indeed, models in which $\sigma$ grows
faster than $\sqrt{\vert x\vert}$ for negative $x$ have typically
infinite bond prices, compare Theorem~4.1 in \cite{Y}. One should
note that all models in Table~\ref{table1} satisfy the conditions of
Assumption~\ref{holder} in the following sense: those models in
which the short rate may fall below zero (V and HW) clearly satisfy
Assumption~\ref{holder}; the remaining
models give rise to non-negative rates, so the parameters can be
chosen arbitrarily for negative values so that
Assumption~\ref{holder} holds without affecting the 
bond value for positive rates.

\section{Auxiliary estimates}
\label{aux}
When dealing with bonds and options on the short rate it is natural
to study the corresponding term structure equation, i.e. the
terminal value problem
\begin{equation}
\label{termstr}
\left\{\begin{array}{ll}
U_t+\alpha U_{xx}+\beta U_x-xU=0 & \mbox{for }t<T\\
U=g& \mbox{for }t=T,\end{array}\right. \end{equation} where
$\alpha:=\sigma^2/2$. Due to technical reasons we instead 
study the function
\begin{equation}
\label{V}
V(x,t):=V^f(x,t):=E_{x,t}\left[\exp\left\{-\int_t^Tf(X_s)\,ds\right\}g(X_T)
\right]
\end{equation}
for some appropriate function $f$ and then let $f$ approach $x$.
Note that the corresponding terminal value problem is
\begin{equation}\label{eqf}
\left\{\begin{array}{ll}
V_t+\alpha V_{xx}+\beta V_x-fV=0 & \mbox{for }t<T\\
V=g & \mbox{for }t=T.\end{array}\right.\end{equation}
Our choice of $f$ is indicated in the following hypothesis.

\begin{hypothesis}
\label{hypf} $f$ is smooth, concave and there exists a constant
$K^\prime>0$ such that
\[f(x)=\left\{\begin{array}{ll}
x & \mbox{if }x\leq K^\prime\\
\mbox{constant} & \mbox{if }x\geq K^\prime+1.\end{array}\right.\]
\end{hypothesis}

The main goal of this section is to provide some estimates
on option prices and the derivatives of $V$, compare
Corollary~\ref{cor} and Proposition~\ref{wdecay}.
We first claim that $V$ grows at most exponentially as $x\to
-\infty$. To be more precise, assume that the pay-off 
function $g$ satisfies
\begin{equation}\label{g}0\leq
g(x)\leq M\max\{ e^{-Kx},1\}
\end{equation}
for some non-negative constants $M$ and $K$ (if $g=1$ then 
\eqref{g} holds
with $M=1$, $K=0$). Moreover, let $D$ be the constant 
from Hypothesis~\ref{hyp}. From \eqref{growthbeta} we have
\begin{equation}\label{below}\beta(x,t)\geq -D(1-x)
\mbox{ for all negative
}x.\end{equation}  
Define
\begin{equation}
\label{W} W(x,t):=e^{f(x)h(t)}V(x,t),
\end{equation}
where
\[h(t)=\frac{e^{D(T-t)}-1}{D}+Ke^{D(T-t)}.\]
The choice of the function $h$ is motivated by the following
observation.
If $\sigma(x,t)=0$ and $\beta(x,t)= Dx$ for negative $x$, and if
$X_0=x\leq 0$, then $X$ satisfies $X_t= xe^{Dt}$. Consequently,
\begin{eqnarray*}
V(x,t) &=& E_{x,t}\left[\exp\left\{-\int_t^Tf(X_s)\,ds\right\}
g(X_T)\right]=\exp\left\{-\int_t^Txe^{D(s-t)}\,ds\right\}
g(xe^{D(T-t)})\\
&\leq& \exp\left\{-x\frac{e^{D(T-t)}-1}{D} \right\}M 
\exp\left\{
-Kxe^{D(T-t)}\right\}=Me^{-xh(t)},
\end{eqnarray*}
where we have used the bound \eqref{g}. This indicates that the
function $W$ is bounded (at least for negative values of
$x$). As is 
shown in Lemma~\ref{growth} below, this intuition is 
indeed correct.

\begin{lemma}
\label{growth} Assume Hypothesis~\ref{hypf} and the bound 
\eqref{g}. Then the function 
$W$ defined in \eqref{W} is bounded on $\R\times [0,T]$, 
i.e. there exists a constant $C$ such that $0\leq W\leq C$.
\end{lemma}

\begin{proof}
By the bound \eqref{g} and the Cauchy-Schwartz inequality we have
\begin{eqnarray}
\label{CS} 
\notag \frac{W^2(x,t)}{M^2} &=& \left(\frac{e^{f(x)h(t)}}{M}
E_{x,t}\left[\exp\left\{-\int_t^Tf(X_s)\,ds\right\}g(X_T)
\right]\right)^2\\
&\leq&
\left(E_{x,t}\left[\exp\left\{ \int_t^T e^{D(s-t)}f(x)-
f(X_s)\,ds\right\} \max\{ e^{-KX_T},1\}e^{Kf(x)e^{D(T-t)}}
\right]\right)^2\\
\notag
&\leq& E_{x,t}\left[\exp\left\{ 2\int_t^T
e^{D(s-t)}f(x)-f(X_s)\,ds\right\}\right]
\\ \notag && \times
E_{x,t}\left[\max\{ e^{-2KX_T},1\}e^{2Kf(x)e^{D(T-t)}} \right].
\end{eqnarray}
By Jensen's inequality (applied to the exponential function and the
time integral) we have
\begin{eqnarray*}
&& \hspace{-15mm}E_{x,t}\left[\exp\left\{ 2\int_t^T
e^{D(s-t)}f(x)-f(X_s)\,ds\right\}\right] \\ &&\leq 
E_{x,t}\left[\frac{1}{T-t}
\int_t^T\exp\left\{2(e^{D(s-t)}f(x)-f(X_s))(T-t)\right\}\,ds\right]
\\ & &= \frac{1}{T-t}
\int_t^TE\left[\exp\left\{-2e^{D(s-t)}(T-t)Y_s\right\}\right] \,ds
\end{eqnarray*}
where $Y_s=e^{-D(s-t)}f(X_s)-f(x)$ with $Y_t=0$. By Ito's 
lemma we have
\[dY_s=\tilde\beta_s\,ds+\tilde\sigma_s\,dB_s \hspace{5mm} 
\mbox{for }s\geq t,\]
where
\[\tilde\beta_s=e^{-D(s-t)}\left(f^\prime(X_s)\beta(X_s,s)
-Df(X_s)+\frac{1}{2} f^{\prime\prime}(X_s,s)\sigma^2(X_s,s)
\right)\] and
\[\tilde\sigma_s= e^{-D(s-t)}f^\prime(X_s)\sigma(X_s,s).\]
It follows from the assumption \eqref{below} and
Hypothesis~\ref{hypf} that $\tilde \beta$ is bounded from below,
i.e. there exists a constant $C$ such that $\tilde \beta_s
\geq -C$ for all $s$ almost surely. Similarly, 
the bound \eqref{growthsigma} and Hypothesis~\ref{hypf}
yield that $\vert\tilde\sigma_s\vert\leq C$. Since $y\mapsto
\exp\{-2e^{D(s-t)}(T-t)y\}$ is decreasing and convex, we have
\begin{equation}
\label{monotone} E\left[\exp\{-2e^{D(s-t)}(T-t)Y_s\}\right]
\leq E\left[\exp\{-2e^{D(s-t)}(T-t)\tilde Y_s\}\right]
\end{equation}
for every $s\in[t,T]$, where
\[\tilde Y_s=- C(s-t)+ CB_{s-t}\]
is a Brownian motion with drift starting from $0$ at time $t$.
Indeed, monotonicity in the drift is immediate, and monotonicity in
the volatility for processes with constant drift and convex pay-off
functions is well-known, compare Theorem~6.2 in \cite{EKJPS} and
Theorem~7 in \cite{JT1}. Consequently,
\begin{eqnarray*}
&&  E_{x,t}\left[\exp\left\{
2\int_t^Tf(x)e^{D(s-t)}-f(X_s)\, ds\right\}\right]\\
&& \hspace{20mm}\leq \frac{1}{T-t} \int_t^TE\left[\exp\left\{-2
e^{D(s-t)}(T-t)\tilde Y_s\right\}\right]\,ds \leq C^\prime
\end{eqnarray*} for some constant $C^\prime$ (in fact, $\tilde Y_s$ is
$N(-C(s-t),C^2(s-t))$-distributed, so explicit bounds of the above
integral are readily derived). Consequently, the first factor in
\eqref{CS} is bounded. For the second factor in \eqref{CS}, note
that
\begin{eqnarray}
\label{secterm}  E_{x,t}\left[\max\{
e^{-2KX_T},1\}e^{2Kf(x)e^{D(T-t)}} \right] &\leq&
E_{x,t}\left[e^{2K(f(x)e^{D(T-t)}-f(X_T))}\right] \\
\notag && + e^{2Kf(x)e^{D(T-t)}}
\end{eqnarray}
The first term in \eqref{secterm} is of the form shown above to be
bounded, and the second one is bounded since $f$ is bounded from
above. This shows that $W$ is bounded, which finishes the proof.
\end{proof}

\begin{corollary}
\label{cor} Assume that $g$ satisfies \eqref{g}. Then the option
price $U$ defined in \eqref{U} above is finite. In fact, 
there exists a constant $M>0$ such that
\[U(x,t)\leq M\max\{ 1, e^{-Mx}\}\] 
for all $(x,t)\in\R\times [0,T]$.
\end{corollary}

\begin{proof}
Since $f(x)\leq x$, we have that $U(x,t)\leq V(x,t)$ for all $x$ and
$t$. Consequently, it follows from Lemma~\ref{growth} that for
every choice of a function $f$ satisfying Hypothesis~\ref{hypf}
there exists a constant $C$ such that $U(x,t)\leq
Ce^{-f(x)(e^{D(T-t)}-1)/D-Ke^{D(T-t)}f(x)}$ for all 
$(x,t)\in\R\times[0,T]$. Thus
\[U(x,t)\leq M\max\{ 1, e^{-Mx}\}\]
for some large constant $M$.
\end{proof}

We next provide conditions under which the $k$th spatial 
derivative of the function $W$ decays like
$\vert x\vert^{-k}$, $k=1,2,3$. We need the following assumptions.

\begin{hypothesis}
\label{hyp} 
The coefficients $\alpha$ and $\beta$ are smooth and $\alpha>0$ at all points. Moreover, 
\begin{equation}
\label{driftcondition} \vert\beta(x,t)-Dx\vert\mbox{ does not 
depend on $x$ for } (x,t)\in\R\times[0,T] \mbox{ with $x$ large negative},
\end{equation}
and the derivatives satisfy the growth conditions
\[\vert\partial_x^k\beta(x,t)\vert\leq C(1+\vert
x\vert)^{1-k}\hspace{12mm}(k=0,1,2,3)\] for $(x,t)\in [0,\infty)\times [0,T]$,
\[\vert\partial_x^k\alpha(x,t)\vert \leq C(1+\vert x\vert)^{2-k}
\hspace{12mm}(k=0,1,2,3)\] for $(x,t)\in[0,\infty)\times[0,T]$, and
\[\vert\partial_x^k\alpha(x,t)\vert
\leq C(1+\vert x\vert)^{-k} \hspace{12mm}(k=0,1,2,3)\] for $(x,t)\in
(-\infty,0]\times [0,T]$.
\end{hypothesis}

\begin{proposition}
\label{wdecay} 
Assume Hypotheses~\ref{hypf} and \ref{hyp}, that the
pay-off function $g$ is smooth, satisfies \eqref{g} and 
that $e^{f(x)K}g(x)$ is constant for large $\vert x\vert$.
Then there exists a
constant $C$ such that the function $W$ defined in \eqref{W}
satisfies
\begin{equation}
\label{hx}
\vert \partial_x^k W \vert\leq \frac{C}{1+\vert x\vert^k}
\end{equation}
for $k=0,1,2,3$.
\end{proposition}

\begin{proof}
The case $k=0$ follows from Lemma~\ref{growth}. Thus we know
that $W$ is a bounded solution to
\begin{equation}
\label{Weq}
\left\{\begin{array}{l}
W_t+\alpha W_{xx}+\hat\beta W_x+\gamma W=0\\
W(x,T)=\hat g(x)\end{array}\right.
\end{equation} 
where
\[\hat g(x)=e^{f(x)K}g(x),\]
\[\hat\beta=\beta-2f_x\alpha h\] and
\[\gamma=(Df-f_x\beta)h+ f_x^2
\alpha h^2- f_{xx}\alpha h.\] 
Let $\hat W(x,t)=W(x,t)-\hat g(x)$. Then
\begin{equation}
\label{What}
\left\{\begin{array}{l}
\hat W_t+\alpha \hat W_{xx}+\hat\beta \hat W_x+\gamma \hat W
+\hat h=0\\
\hat W(x,T)=0,\end{array}\right.
\end{equation} 
where 
\[\hat h=\alpha \hat g_{xx}+\hat\beta \hat g_x+\gamma \hat g.\]
Since $\hat g$ is constant for large values of $\vert x\vert$, 
$W$ satisfies \eqref{hx} if and only if $\hat W$ does.
Using Hypotheses~\ref{hypf} and \ref{hyp} we note that
\[\begin{array}{ll}
\vert\partial_x^k\alpha\vert \leq C(1+\vert x\vert)^{2-k} &
\vert\partial_x^k\hat\beta\vert\leq C(1+\vert
x\vert)^{1-k}\\
\vert\partial_x^k\gamma\vert\leq C(1+\vert x\vert)^{-k}
& \vert\partial_x^k \hat h\vert\leq C(1+\vert x\vert)^{-k}
\end{array}\] for some constant $C$ and for all
$k=0,1,2,3$. By the Feynman-Kac representation theorem,
\[\hat W(x,t)=E\left[\int_t^T\exp\left\{\int_t^s\gamma(Y^{x,t}_r,r)\,dr
\right\} \hat h(Y^{x,t}_s)\,ds\right],\] where $Y^{x,t}$ is the diffusion
\[dY_s^{x,t}=\hat\beta(Y^{x,t}_s,s)\,ds+\sigma(Y^{x,t}_s,s)\,dB_s\]
with $Y^{x,t}_t=x$. Now, let $x<y$. Then
\begin{eqnarray*}
\vert \hat W(x,t)-\hat W(y,t) \vert &\leq& 
E\left[\int_t^T\exp\left\{\int_t^s\gamma(Y^{x,t}_r,r)\,dr
\right\}\left \vert \hat h(Y^{x,t}_s)-\hat h(Y^{y,t}_s)\right \vert\,ds\right]\\
&&\hspace{-20mm}+ E\left[\int_t^T\left\vert\exp\left\{\int_t^s\gamma(Y_r^{x,t},r)\,dr\right\}
-\exp\left\{\int_t^s\gamma(Y_r^{y,t},r)\,dr\right\}\right\vert 
\hat h(Y^{y,t}_s)\,ds\right]\\
&=& I_1+I_2.
\end{eqnarray*}
Since $\gamma$ is bounded and $\hat h$ is Lipschitz continuous, 
we find that
\[I_1 \leq C\int_t^T E\left[\vert Y_s^{x,t}-Y_s^{y,t}\vert\right]\,ds=
C\int_t^T E\left[Y_s^{y,t}-Y_s^{x,t}\right]\,ds,\]
where the equality holds since $x<y$ implies $Y^{x,t}_s\leq
Y^{y,t}_s$ for $s\geq t$, compare Theorem~IX.3.7 in \cite{RY}. 
Since the drift $\hat\beta$ is Lipschitz continuous, we have
\[E\left[Y_s^{y,t}-Y_s^{x,t}\right]=y-x+\int_t^s E\left[
\beta(Y_s^{y,t},r)-\beta(Y_r^{x,t},r)\right]\,dr
\leq y-x+C\int_t^s E\left[Y^{y,t}_r-Y^{x,t}_r\right]\,dr\]
for any $s\geq t$. It follows from Gronwall's lemma that 
\[I_1\leq C\vert y-x\vert(T-t).\]
Similarly, since $\hat h$ is bounded and $\gamma$ is bounded and 
Lipschitz it can also be shown that
\[I_2\leq C\vert y-x\vert(T-t)\]
for some constant $C$. It follows that $\hat W_x$ is bounded on 
$\R\times [0,T]$ and that 
\[\lim_{(x,t)\to(x_0,T)}\hat W_x(x,t)=0.\]
Therefore, by differentiating the equation \eqref{What} we know that the first spatial derivative $p:=\hat W_x$ is a bounded 
classical solution to
\begin{equation}
\label{p}
\left\{\begin{array}{l}
p_t+\alpha p_{xx}+(\alpha_x+\hat \beta) p_x+(\gamma+
\hat \beta_x)p+\gamma_x \hat W-\hat h_x=0\\
p(x,T)=0\end{array}\right.
\end{equation}
It is straightforward to check that
\[\tilde p(x,t):=e^{C(T-t)}l(x),\] where $l$ is a smooth positive
function with $l(x)=1/(1+\vert x\vert)$ for large $\vert x\vert$, is
a supersolution to \eqref{p} for some large constant $C$, and
$-\tilde p$ is a subsolution. By the maximum principle we have $\vert \hat W_x\vert=\vert p\vert\leq \tilde p \leq
e^{CT}/(1+\vert x\vert)$. This proves \eqref{hx}
for $k=1$.

Next, the solution to equation \eqref{p} also has a stochastic representation. Using this representation, the same analysis
as above can be applied to show that $p_x$ is bounded, and 
using the maximum principle again it is straightforward
to check that $p_x=\hat W_{xx}$ decays like $\vert x\vert^{-2}$. Finally, 
the same reasoning can be applied to prove \eqref{hx} for $k=3$.
\end{proof}

\section{Continuity of bond prices in the model parameters}
\label{continuity}

As mentioned in the introduction, it is central to our
approach to be able to approximate a bond price (or an 
option price) by a sequence of bond prices (option prices) 
in models with regular coefficients. Thus we need to know 
that bond prices are continuous in the model parameters. 
To formulate such a result, let $\sigma_n$ and $\beta_n$ satisfy 
the Assumption~\ref{holder} uniformly in $n$ (i.e.
the bounds \eqref{growthsigma} and \eqref{growthbeta} hold 
with the
same constant $D$). Also assume that $\alpha_n\to \alpha$ and
$\beta_n\to\beta$ uniformly on compacts as $n\to\infty$. 
Let $X^n$ be the solution of the stochastic 
differential equation
\[dX^n_s=\beta_n(X^n_s,s)\,ds + \sigma_n(X^n_s,s)\,dB_s
\hspace{10mm} X^n_t=x.\]
Under these conditions it is known that
\begin{equation}
\label{sup} E_{x,t}\left[\sup_{t\leq s\leq T}\vert
X^n_s-X_s\vert^2\right]\to 0 \end{equation} as $n\to\infty$,
see Theorem~2.5 in \cite{BMO}. Now let
\begin{equation}
\label{V_n} U_n(x,t)=E_{x,t}\left[\exp\left\{-\int_t^T
X^n_s\,ds\right\} g(X_T^n)\right].\end{equation} Here $g$ is 
assumed to have the following property: for each positive 
constant $k$ there exists a constant $C_k$ such that
\begin{equation} \label{ql}\vert g(x)\wedge k-
g(y)\wedge k\vert\leq
C_k\vert x-y\vert\end{equation} 
for all $x,y\in \R$. Note that this condition is satisfied 
for example if $g$ is convex and decreasing.

\begin{proposition}
\label{approx} Assume that $\sigma_n$, $\beta_n$ and $X^n$ are as
described above. Also assume that $g$ satisfies \eqref{g} and
\eqref{ql}. Then
\[U(x,t)=\lim_{n\to\infty}U_n(x,t)\]
for all $(x,t)\in\R\times[0,T]$.
\end{proposition}

\begin{proof}
Fix a point $(x,t)\in \R\times [0,T]$. For positive constants $k$,
define $U^k(x,t)$ and $U_n^k(x,t)$ by
\[U^k(x,t)=E_{x,t}\left[\exp\left\{\int_t^T k\wedge
-X_s\,ds\right\} (g(X_T)\wedge k) \right]\] and
\[U^k_n(x,t)=E_{x,t}\left[\exp\left\{\int_t^T k\wedge -X^n_s\,ds\right\}
(g(X^n_T)\wedge k) \right].\] Then
\[\left\vert U_n^k(x,t)-U^k(x,t)\right\vert \leq I_1+I_2,\]
where
\[I_1 = E_{x,t}\left[ \exp\left\{\int_t^Tk\wedge
-X^n_s\,ds\right\}\left\vert g(X_T^n)\wedge k-g(X_T)\wedge
k\right\vert\right]\] and 
\[I_2 = E_{x,t} \left[\left\vert \exp\left\{\int_t^Tk\wedge
-X_s^n\,dt\right\}-\exp\left\{\int_t^Tk\wedge -X_s\,
ds\right\}\right\vert (g(X_T)\wedge k)\right].\] Note that
\[I_1\leq e^{k(T-t)}E_{x,t}\left[\left\vert g(X_T^n)\wedge
k-g(X_T)\wedge k\right\vert\right]\leq e^{k(T-t)}
C_kE_{x,t}\left[\left\vert X_T^n-X_T\right\vert\right]\] 
since $g$ satisfies \eqref{ql}. It thus follows from 
\eqref{sup} that $I_1\to 0$ as $n\to\infty$. 
Similarly, using the inequality $\vert
e^x-e^y\vert\leq (e^x+e^y)\vert x-y\vert$, we find
\begin{eqnarray*}
I_2 &\leq& kE_{x,t}\left[ \left\vert
\exp\left\{\int_t^Tk\wedge
-X_s^n\,ds\right\}-\exp\left\{\int_t^Tk\wedge -X_s\,
ds\right\}\right\vert\right]\\
&\leq& 2ke^{k(T-t)}E_{x,t}\left[\int_t^T\vert k\wedge
-X_s^n-k\wedge -X_s\vert\,ds\right]\\
&\leq& 2ke^{k(T-t)}E_{x,t}\left[(T-t)\sup_{t\leq s\leq T}\vert
X_s^n- X_s\vert\right]\to 0
\end{eqnarray*}
as $n\to\infty$ by \eqref{sup}. Consequently,
\begin{equation}
\label{lim} \lim_{n\to\infty}U_n^k (x,t) = U^k (x,t)
\end{equation}
for each fixed $k$. It then follows from the monotone convergence
theorem that
\begin{eqnarray}
\label{lim2} \lim_{k\to\infty}\lim_{n\to\infty}U^k_n(x,t) &=&
\lim_{k\to\infty} E_{x,t}\left[\exp\left\{\int_t^T k\wedge
-f(X_s)\,ds\right\}
(g(X_T)\wedge k)\right]\\
\notag &=& E_{x,t}\left[\exp\left\{-\int_t^T f(X_s)\,ds\right\}
g(X_T) \right]=U(x,t).
\end{eqnarray}
Since
\[U_n(x,t)=\lim_{k\to\infty}U^k_n(x,t)\]
by the monotone convergence theorem, it suffices to show that the
order of the limits in \eqref{lim2} can be interchanged. To see
this, first note that
\begin{eqnarray}
\label{est3} 0 &\leq& U_n(x,t)-U^k_n(x,t)\\ \notag &\leq&
E_{x,t}\left[\exp\left\{-\int_t^T X^n_s\,ds\right\}g(X_T^n)
-\exp\left\{\int_t^T k\wedge -X^n_s\,ds\right\}g(X_T^n)
\right]\\
\notag && +E_{x,t}\left[\exp\left\{-\int_t^T X_s^n\,ds
\right\}(g(X_T^n)-k\wedge g(X_T^n))\right]=I_3+I_4.\end{eqnarray} We
first claim that
\begin{equation}
\label{I3} \lim_{k\to\infty}\sup_{n} I_i=0\hspace{10mm}\mbox{for
}i=3,4.
\end{equation}
To see this, note that
\begin{eqnarray*}
I_3 &=& E_{x,t}\left[\exp\left\{-\int_t^T X^n_s \,ds\right\}
g(X_T^n)\left(1-\exp
\left\{\int_t^T\min(k+X_s^n,0)\,ds\right\}\right) \right]\\
\notag &\leq& E_{x,t}\left[\exp\left\{-\int_t^T X^n_s\,ds\right\}
g(X_T^n)
\int_t^T (-k-X_s^n)^+\,ds\right]\\
\notag &\leq& \left(E_{x,t}\left[\exp\left\{-2\int_t^T
X^n_s\,ds\right\}(g(X_T^n))^2\right] 
E_{x,t}\left[ \left(\int_t^T(-k-X_s^n)^+\,ds\right)^2
\right]\right)^{1/2}
\end{eqnarray*}
where we have used the inequality $1-e^x\leq -x$ and the
Cauchy-Schwartz inequality. The first factor is bounded 
uniformly in $n$ for each fixed $x$ (this can be proved 
analogously to Lemma~\ref{growth}). As for the second factor, 
an application of Jensen's inequality gives
\begin{equation*}
E_{x,t}\left[ \left(\int_t^T(-k-X_s^n)^+\,ds\right)^2 \right]\leq
(T-t)\int_t^T
E_{x,t}\left[\left((-k-X^n_s)^+\right)^2\right]\,ds.\end{equation*}
Proceeding as in the proof of Lemma~\ref{growth} above, it is
straightforward to check that
\begin{eqnarray*}
\label{est} E_{x,t}\left[((-k-X^n_s)^+)^2\right] &\leq&
e^{2DT}E_{x,t}\left[\left((-ke^{-DT}-e^{-D(s-t)}X_s^n)^+\right)^2\right]
\\
\notag &\leq& e^{2DT}E_{x,t} \left[\left((-ke^{-DT}-\tilde
Y_s)^+\right)^2\right]
\end{eqnarray*}
for a Brownian motion $\tilde Y$ with negative drift. Moreover,
since $((-ke^{-DT}-\tilde Y_s)^+)^2$ is a submartingale we have
\begin{equation*}
\label{unif} E_{x,t}\left[\left((-ke^{-DT}-\tilde
Y_s)^+\right)^2\right]\leq E_{x,t}\left[\left((-ke^{-DT}-\tilde
Y_T)^+\right)^2\right]\to 0\end{equation*} as $k\to\infty$ by
dominated convergence. Consequently, \eqref{I3} holds for $i=3$.
Similar methods can be employed to establish \eqref{I3} also for
$i=4$.

Let $\ep>0$ be given. From \eqref{est3} and \eqref{I3} it follows
that there exists a $k_0$ such that
\[0\leq U_n(x,t)-U_n^k(x,t)\leq \ep\]
for all $k\geq k_0$ and all $n$. Thus, for such a $k$ we have
\[\lim_{n\to\infty} U^k_n(x,t)\leq
\liminf_{n\to\infty} U_n(x,t) \leq \limsup_{n\to\infty} U_n(x,t)
\leq \lim_{n\to\infty} U_n^k(x,t)+\ep\] where we recall from
\eqref{lim} that the outer limits exist. Letting $k\to\infty$ and
remembering \eqref{lim2} above we find that
\[U(x,t)\leq \liminf_{n\to\infty}U_n(x,t) \leq
\limsup_{n\to\infty}U_n(x,t)\leq U(x,t)+\ep.\] Since $\ep$ is
arbitrary,
\[\lim_{n\to\infty}U_n(x,t)=U(x,t)\]
as required.
\end{proof}

\section{Convexity of bond prices}
\label{convexity} 
In this section we show that interest rate models
have convex bond prices provided the drift $\beta$ is not ``too
convex''. More precisely, we need to require that $\beta_{xx}\leq 2$
(in the sense of distributions). To see why this condition comes
into play, consider the corresponding term structure equation, i.e.
the terminal value problem
\[\left\{\begin{array}{l}
U_t+\alpha U_{xx}+\beta U_x-xU=0\\
U(x,T)=g(x)\end{array}\right.\] for some convex and 
decreasing pay-off function $g$. When using the PDE-approach we will
in the sequel simplify the presentation by performing a standard
change $\tau= T-t$ of the direction of time. By a slight abuse of
notation, we use the same symbols $\beta$, $\sigma$, $U$, $V$ etc.
to denote the new functions depending on time $\tau$ to 
maturity 
rather than on the actual time $t$. With this new convention, the
term structure equation becomes an {\em initial} value problem
\[\left\{\begin{array}{l}
U_\tau=\alpha U_{xx}+\beta U_x-xU\\
U(x,0)=g(x).\end{array}\right.\] Assume for the moment that all
coefficients are regular enough. Also assume that there is
a first point $(x_0,\tau_0)$ at which convexity is almost 
lost, i.e. that
$U(x,\tau)$ is convex for $0\leq \tau\leq \tau_0$ and
$U_{xx}(x_0,\tau_0)=0$. Then
\begin{eqnarray*}
\partial_\tau U_{xx} &=& \partial_x^2 U_\tau=
\partial_x^2 (\alpha U_{xx}+\beta U_x-xU)\\
&=& \alpha U_{xxxx}+(2\alpha_x+\beta)U_{xxx}+(2\beta_x-x)U_{xx}
+(\beta_{xx}-2)U_x.
\end{eqnarray*}
Since $x\mapsto U_{xx}(x,\tau_0)$ has a minimum at $x=x_0$, we have
$U_{xxxx}\geq U_{xxx}=U_{xx}=0$ at $(x_0,\tau_0)$. Thus we find that
\[\partial_\tau U_{xx}\geq (\beta_{xx}-2)U_x\geq 0\] at this point
provided $\beta_{xx}\leq 2$, i.e. the infinitesimal change of $U$ at
the critical point $(x_0,\tau_0)$ is convex. This suggests that
convexity will be preserved if $\beta_{xx}\leq 2$. Below we make
this argument rigorous.

\begin{theorem}
\label{main} Assume that $\beta_{xx}(x,t)\leq 2$ (in the sense of
distributions) at all points $(x,t)\in\R\times [0,T]$. Also assume
that the pay-off function $g$ is convex, decreasing and satisfies
\eqref{g}. Then the corresponding option price $U(x,t)$ is convex in
$x$ at all times $t\in[0,T]$.
\end{theorem}

\begin{remark}
Note that all models in Table~\ref{table1} satisfy the condition
$\beta_{xx}\leq 2$. Consequently, all those models give rise to
convex bond prices.
\end{remark}

\begin{proof}
We first assume that Hypothesis~\ref{hyp} holds, that $\beta_{xx}=0$
for $x\geq C$ for some constant $C$, and that the
pay-off function $g$ is smooth with $e^{Kf(x)}g(x)$ being constant
for large $\vert x\vert$, where $K$ is the constant from 
\eqref{g}.

Instead of studying the option price $U$ directly, we first study
the function $V=V^f$ defined in \eqref{V} for some function $f$
satisfying Hypothesis~\ref{hypf} such that
\begin{equation}
\label{betaxx} \beta_{xx}(x,\tau)-2f_x(x)\leq 0
\end{equation}
for all $(x,\tau)\in \R\times [0,T]$. Under these assumptions we
claim that the function $V$ is spatially convex.

To see this, note that it follows from Proposition~\ref{wdecay} 
that there exists a constant $C$ such that
\begin{equation}
\label{Vxx}
\vert V_{xx}(x,\tau)\vert\leq Cp(x,\tau)
\end{equation}
at all points, where \[p(x,\tau)=e^{-f(x)K_0e^{(D+1)\tau}}\] 
(here $K_0$ is a large constant so that 
$K_0\geq (e^{DT}-1)/D+Ke^{DT}$). For $\ep>0$,
consider the function
\[V^\ep(x,\tau):=V(x,\tau)+\ep e^{M\tau}(x^2+x^Np(x,\tau)).\]
Here the even number $N>2$ is chosen large so that
$x^2+x^Np(x,\tau)$ has a strictly positive second spatial
derivative at all points $(x,\tau)\in\R\times [0,T]$ 
(the constant
$M$ will be chosen below). Since $V$ satisfies $\partial_\tau V
=\mathcal L^f V$ where
\[\mathcal L^f=\alpha\partial_x^2+\beta\partial_x-f,\] we find that
\begin{equation}
\label{eq1}  \partial_x^2(\partial_\tau V^\ep-\mathcal L^f V^\ep) =
\ep \partial_x^2(\partial_\tau -\mathcal L^f)
\left(e^{M\tau}(x^2+x^N p)\right) =\ep e^{M\tau}(MI_1-I_2),
\end{equation}
where
\[I_1=\partial_x^2(x^2+x^Np)>0\]
and
\[I_2=\partial_x^2\left(\alpha\partial_x^2(x^2+x^Np)+\beta
\partial_x(x^2+x^Np)-f(x^2+x^Np)+(D+1)K_0e^{(D+1)\tau}
fx^N\right).\]
For large positive $x$ we have that $f$ is constant, so $p$ is
bounded. Thus the term $I_1$ grows like $x^{N-2}$, whereas the term
$I_2$ grows at most like $x^{N-2}$ for large positive $x$.
Similarly, for large negative $x$ we have $f(x)=x$, so
\[I_1 \sim K_0^2e^{2(D+1)\tau}x^Np, \] 
whereas 
\[I_2 \sim \alpha x^N p_{xxxx}+
(-\beta K_0 e^{(D+1)\tau}-x+(D+1)K_0e^{(D+1)\tau}x)x^Np_{xx}.\]
Using \eqref{growthbeta} we find that the highest order
terms of $I_2$ behaves at most like
\[\alpha x^N p_{xxxx}+
(-D K_0 e^{(D+1)\tau}-1+(D+1)K_0e^{(D+1)\tau})x^{N+1}p_{xx},\]
and since $-D K_0 e^{(D+1)\tau}-1+(D+1)K_0e^{(D+1)\tau}\geq 0$
we find that there exists a positive constant $C$ such that
\[I_2\leq Cx^Np\]
for large negative $x$. Consequently, $I_1$ grows at least
as fast as $I_2$ for large values of $\vert x\vert$,
so $M$ can be chosen large so that $MI_1$ dominates $I_2$ 
everywhere. Actually, we will assume that $M$ is chosen so large that
\begin{equation}
\label{in} M I_1>I_2-I_3
\end{equation}
at all points, where
\[I_3:=(\beta_{xx}-2f_x)\partial_x(x^2+x^Np)\]
(compare \eqref{spaceder} below). This can be done since $\beta_{xx}=0$ and thus $I_3=0$ for large $x$, and
$I_3\geq 0$ for large negative $x$, compare \eqref{betaxx}.

Now, let
\[\Gamma:=\{(x,\tau):V^\ep_{xx}(x,\tau)<0\}.\]
We claim that the set $\Gamma$ is empty. To see this, note 
that the second spatial derivative of $x^2+x^Np$ 
grows at least like
$Cx^Np$ as $\vert x\vert\to\infty$. Consequently, it follows
from \eqref{Vxx} that $\Gamma\subseteq [-R,R]\times [0,T]$ for some $R>0$. 
Thus $\Gamma$
is bounded and $\overline \Gamma$ is compact. Suppose that
$\Gamma\not=\emptyset$, and define
\[\tau_0:=\inf\{\tau\geq 0:(x,\tau)\in\overline \Gamma
\mbox{ for some }x\in\R\}.\]
Due to compactness, the infimum is attained at $(x_0,\tau_0)$ for some $x_0$, and by continuity of $V_{xx}$ we have $V^\ep_{xx}(x_0,\tau_0)=0$, so $\tau_0>0$.
Since $V^\ep_{xx}(x_0,\tau)\geq 0$ for $\tau\leq \tau_0$, we must have
\begin{equation}
\label{timeder}
\partial_x^2\partial_\tau V^\ep(x_0,\tau_0)= \partial_\tau
\partial_x^2 V^\ep(x_0,\tau_0)
\leq 0.\end{equation} Moreover, straightforward calculations give
\begin{eqnarray*}
\partial_x^2 (\mathcal L^f V^\ep) &=&\partial_x^2 (\alpha V^\ep_{xx}+\beta
V^\ep_x -fV^\ep)\\
&=& \alpha V^\ep_{xxxx}+(2\alpha_x+\beta)V^\ep_{xxx}+
(\alpha_{xx}+2\beta_x-f)V^\ep_{xx}\\
&& +(\beta_{xx}-2f_{x})V^\ep_x -f_{xx}V^\ep.
\end{eqnarray*}
Since $V^\ep_{xx}(x_0,t_0)=0$ and $V^\ep$ is convex, the function
$x\mapsto V^\ep_{xx}(x,\tau_0)$ has a minimum point at $x_0$.
Consequently, $V^\ep_{xxx}(x_0,\tau_0)=0$ and
$V^\ep_{xxxx}(x_0,\tau_0)\geq 0$. Therefore, at the point
$(x_0,\tau_0)$ we have
\begin{eqnarray}
\label{spaceder}
\partial_x^2 (\mathcal L^f V^\ep) &=& \alpha
V^\ep_{xxxx}+(\beta_{xx}-2f_x)V^\ep_x -f_{xx}V^\ep \\
\notag &\geq&  \ep e^{M\tau}(\beta_{xx}-2f_x)
\partial_x(x^2+x^Np)=\ep e^{M\tau}I_3
\end{eqnarray}
where we have used $V_x\leq 0$, \eqref{betaxx} and that $f$ is
concave. Combining \eqref{eq1}, \eqref{in}, \eqref{timeder} and
\eqref{spaceder} yields
\[\partial_x^2(\partial_\tau V^\ep-\mathcal L^f V^\ep)=\ep e^{M\tau}
(MI_1-I_2)>-\ep e^{M\tau} I_3\geq
\partial_x^2(\partial_\tau V^\ep-\mathcal L^f V^\ep)\]
at $(x_0,\tau_0)$. This contradiction shows that the set $E$ is
empty, i.e. $V^\ep$ is convex. Letting $\ep\downarrow 0$ it follows
that also $V$ is spatially convex at all times $\tau\in[0,T]$.

To deduce the convexity of $U$, consider an increasing sequence
$\{f_i\}^\infty_{i=1}$ of functions $f_i$ satisfying
Hypothesis~\ref{hypf} such that $f_i(x)\to x$ as $i\to\infty$ for
all $x$. Then
\[\int_t^T f_i(X_s)\,ds\to \int_t^TX_s\,ds\]
almost surely as $i\to\infty$. Therefore
\begin{eqnarray*}
V^{f_i}(x,t) &=& E_{x,t}\left[\exp\left\{
-\int_t^Tf_i(X_s)\,ds\right\} g(X_T)\right] \\
&\to& E_{x,t}\left[\exp\left\{
-\int_t^TX_s\,ds\right\}g(X_T)\right]=U(x,t)
\end{eqnarray*} 
as $i\to\infty$  by
monotone convergence. Since each function $V^{f_i}$ is spatially
convex it follows that the option price $U(x,t)$ is spatially
convex.

To remove the Hypothesis~\ref{hyp} and the assumption that
$\beta_{xx}=0$ for large $x$, assume only Assumption~\ref{holder} and
that $\beta_{xx}\leq 2$ (in the sense of distributions). Let
$\{\sigma_n\}_{n=1}^\infty$ and $\{\beta_n\}_{n=1}^\infty$ be
sequences of smooth continuous coefficients satisfying Hypothesis~\ref{hyp}
such that $(\beta_n)_{xx}=0$ for large $x$ and $(\beta_n)_{xx}\leq 2$.
Moreover, assume that $\sigma_n$ and $\beta_n$ satisfy the growth
conditions \eqref{growthsigma} and \eqref{growthbeta} uniformly in
$n$ and that $\sigma_n\to\sigma$ and $\beta_n\to\beta$ uniformly on
compacts as $n\to\infty$. It follows from Proposition~\ref{approx}
that
\[U(x,t)=\lim_{n\to\infty}U_n(x,t),\]
where the function $U_n$ is defined as in \eqref{V_n}. Since the
pointwise limit of a sequence of convex functions is convex, it
follows that $U$ is convex.

Finally, to remove the assumptions about the smoothness of $g$ and
that $e^{Kf(x)}g(x)$ is constant outside a compact we approximate $g$ from above by a sequence $\{g_n\}_{n=1}^\infty$ of smooth pay-offs behaving like $e^{-Kf(x)}$ outside compacts, such that $g_n(x)\downarrow
g(x)$ as $n\to\infty$. By monotone convergence it follows
that $U$ is convex also without the smoothness requirements on
$g$.
\end{proof}

The heuristic calculations in the beginning of this section indicate
that the condition $\beta_{xx}\leq 2$ is not only a sufficient
condition, but also a {\em necessary} condition for preservation of
convexity for the term structure equation. Our next result shows
that this is indeed true provided the coefficients are regular
enough.

\begin{theorem}
\label{nec} Assume that $\alpha$ and $\beta$ are smooth. Also assume that $\alpha>0$ and $\beta_{xx}>2$ at some point $(x_0,T)$. Then there exists an option
with maturity $T$ and with a decreasing and convex pay-off $g$ such
that the corresponding price $U(x,t)$ is non-convex at some time
$t<T$.
\end{theorem}

\begin{proof}
Let $g$ be a smooth convex pay-off function which is linear and
strictly decreasing in a neighborhood of $x_0$. Since 
$U(x,\tau)$ is a
solution of a parabolic differential equation with regular
coefficients, its derivatives exist and are continuous up to the
boundary $\tau=0$, see \cite{L} (here we again let $\tau=T-t$). Straightforward calculations yield
\begin{eqnarray*}
\partial_\tau U_{xx} &=& \partial_x^2 U_\tau
=\partial_x^2(\alpha U_{xx}+\beta U_x -x U)\\
&=& \alpha U_{xxxx}+(2\alpha_x+\beta)U_{xxx}+(\alpha_{xx}+2\beta_x)
U_{xx} + (\beta_{xx}-2)U_x.
\end{eqnarray*}
Since $g$ is linear in a neighborhood of $x_0$, we find that
\[\partial_\tau U_{xx}(x_0,0)=(\beta_{xx}-2)g_x>0.\]
Since $U_{xx}(x_0,0)=0$, this means that $U$ is not convex at some
time $\tau>0$. This finishes the proof.
\end{proof}

\section{Parameter monotonicity}
\label{mon}

In this section we utilize the well-known connection between
convexity and parameter monotonicity, see for instance
\cite{EKJPS} or \cite{ET1}. We thus show how preservation 
of convexity implies that
bond prices and prices of convex options are monotonically
increasing in the volatility. To formulate the result, let $X$ and
$\tilde X$ be two diffusion processes satisfying
\[dX_t=\beta(X_t,t)\,dt+\sigma(X_t,t)\,dB_t\]
and
\[d\tilde X_t=\tilde\beta(\tilde X_t,t)\,dt+\tilde\sigma(\tilde
X_t,t)\,dB_t,\] respectively. Let
\[U(x,t)=E_{x,t}\left[\exp\left\{-\int_t^T X_s\,ds\right\}g(X_T)\right]\]
and
\[\tilde U(x,t)=E_{x,t} \left[\exp\left\{-\int_t^T
\tilde X_s\,ds\right\}g(\tilde{X}_T)\right]\] be the corresponding
option prices.

\begin{theorem}
\label{mono} Assume that $\beta(x,t)\leq \tilde \beta(x,t)$ and
$\vert\sigma(x,t)\vert\geq \vert\tilde\sigma(x,t)\vert$ for all
$(x,t)\in\R\times [0,T]$, and that either $\beta_{xx}\leq 2$ at all
points or $\tilde\beta_{xx}\leq 2$ at all points (both in the
distributional sense). Also assume that the pay-off function $g$ is
convex, decreasing and satisfies \eqref{g}. Then $\tilde U(x,t)\leq
U(x,t)$ for all $(x,t)\in\R\times [0,T]$.
\end{theorem}

\begin{remark}
Monotonicity in the drift is immediate since $g$ is decreasing and
it is known that $\beta(x,t)\leq \tilde \beta(x,t)$ for all $x$ and
$t$ implies that $X_T\leq \tilde X_T$ (if
$\sigma(x,t)=\tilde\sigma(x,t)$ for all $x$ and $t$).
Theorem~\ref{mono} tells us that bond prices are also increasing in
the volatility, provided the condition $\beta_{xx}\leq 2$ 
is fulfilled. In particular, all models in Table~\ref{table1}
are monotonically increasing in the volatility and decreasing in the drift.
\end{remark}

\begin{proof}
Let $f$ be a given function satisfying Hypothesis~\ref{hypf}, and define the functions
\[V(x,t)=E_{x,t}\left[\exp\left\{-\int_t^Tf(X_s)\,ds\right\}g(X_T)
\right]\] and
\[\tilde V(x,t)=E_{x,t}\left[\exp\left\{-\int_t^Tf(\tilde X_s)\,ds\right\}
g(\tilde X_T) \right].\] As in the proof of Theorem~\ref{main},
we will use the parametrization in remaining time
$\tau=T-t$ to maturity. Let
\[p(x,\tau)=e^{-f(x)K_0e^{(D+1)\tau}},\]
where the positive constant $D$ satisfies \eqref{growthsigma}
and \eqref{growthbeta} for both models, and the constant
$K_0$ is large so that $K_0\geq
(e^{DT}-1)/D+Ke^{DT}$. For $\ep>0$, define
\[V^\ep(x,\tau)=V(x,\tau)+\ep e^{M\tau}(1+x^2p(x,\tau)).\]
It is straightforward to check that the constant $M$ can be
chosen so large that
\begin{equation}
\label{estimate} M(1+x^2p) > -x^2\partial_\tau p+ \alpha\partial_x^2(1+x^2p)+
\beta\partial_x(1+x^2p)-f(1+x^2p)
\end{equation} at all points $(x,\tau)\in\R\times
[0,T]$. Now, let
\[\Gamma:=\{(x,\tau)\in\R\times [0,T]:\tilde V(x,\tau)<V^\ep(x,\tau)\},\]
and suppose that $\Gamma\not=\emptyset$. Since 
$\tilde V p^{-1}$ is
bounded by Lemma~\ref{growth}, and since $V^\ep p^{-1}$ grows
like $x^2$ for large $\vert x\vert$, there exists $R>0$ such that
$\Gamma\subseteq (-R,R)\times [0,T]$. Therefore, $\Gamma$ is 
bounded and $\overline \Gamma$ is compact. Define
\[\tau_0:=\inf\{\tau:(x,\tau)\in\overline \Gamma 
\mbox{ for some }x\}.\]
By compactness, there exists $x_0\in\R$ such that
$(x_0,\tau_0)\in\overline \Gamma$. Since $V^\ep-\tilde V$ is continuous,
we have $V^\ep(x_0,\tau_0)=\tilde V(x_0,\tau_0)$. Therefore,
$V^\ep(x,0)>\tilde V(x,0)$ implies that $\tau_0>0$. By the
definition of $\tau_0$ we have $V^\ep(x_0,\tau)-\tilde V(x_0,\tau)\geq 0$ for $0<\tau<\tau_0$, so
\begin{equation}
\label{time}
\partial_\tau (V^\ep-\tilde V)\leq 0
\end{equation} 
at $(x_0,\tau_0)$.
On the other hand, $V$ and $\tilde V$ satisfy the parabolic
equations $V_\tau=\mathcal L^f V$ and 
$\tilde V_\tau=\tilde{\mathcal L}^f\tilde V$, respectively, 
where
\[\mathcal L^f=\alpha \partial_x^2+\beta \partial_x-f\hspace{5mm}\mbox{
and}\hspace{5mm} \tilde{\mathcal L}=\tilde\alpha \partial_x^2+
\tilde\beta
\partial_x-f.\]
Consequently, at the point $(x_0,\tau_0)$ we have
\begin{eqnarray*}
\partial_\tau (V^\ep-\tilde V) &=& \mathcal L^fV^\ep-
\tilde{\mathcal L}^f \tilde V +\ep e^{M\tau}M(1+x^2p)\\
&&+\ep e^{M\tau}\left(x^2\partial_\tau p-\alpha \partial_x^2(1+x^2p)
-\beta\partial_x(1+x^2p)+f(1+x^2p) \right)\\ 
&>& \alpha V^\ep_{xx} +\beta V^\ep_x-
fV^\ep-(\tilde\alpha \tilde V_{xx} +\tilde\beta\tilde V_x -f\tilde
V)
\end{eqnarray*}
by \eqref{estimate}. The function $x\mapsto V^\ep(x,\tau_0)-\tilde
V(x,\tau_0)$ attains its minimum 0 at $x=x_0$. Thus we have
$V^\ep=\tilde V$, $V^\ep_x=\tilde V_x$ and $V^\ep_{xx}\geq \tilde
V_{xx}$ at the point $(x_0,\tau_0)$. Since at least one of the two
models is convexity preserving, $V^\ep_{xx}\geq 0$ at this point, so
$\alpha V^\ep_{xx}\geq \tilde\alpha\tilde V_{xx}$. Therefore, at
$(x_0,\tau_0)$,
\[\partial_\tau (V^\ep-\tilde V)>\alpha V^\ep_{xx}-\tilde\alpha \tilde
V_{xx}+ (\beta-\tilde\beta)\tilde V_x\geq 0,\] 
where we also have used
$\beta\leq \tilde\beta$ and $\tilde V_x\leq 0$. This contradicts
\eqref{time}. Consequently, $\Gamma$ is empty, so $V^\ep\geq \tilde V$
everywhere. Letting $\ep\searrow 0$ it follows that $V\geq \tilde
V$.

Finally, an application of the monotone convergence theorem (as $f$
approaches $x$) yields that $U\geq \tilde U$.
\end{proof}

\section{Bond options}
\label{bond options}

In this section we apply the results of Sections~\ref{convexity} and
\ref{mon} to study convexity and monotonicity properties of bond
option prices. We consider a European call option with time of
maturity $T_1$ and strike price $K>0$ on a bond maturing at $T_2$,
where naturally $T_2>T_1$. We denote the price of this option at
time $t\le T_1$ by $C(x,t; T_1, T_2)$ with $x$ denoting the short
rate. Thus
\[C(x,t;T_1,T_2)=E_{x,t}\left[\exp\left\{-\int_t^{T_1}
X_s\,ds\right\}(u(X_{T_1},T_1)-K)^+\right],\] where $U(x,t)$ is the
value function of a $T_2$-bond. Assuming that the short rate
dynamics is given by equation (\ref{X}), we then have the following
result.

\begin{theorem}
\label{callconvexity}  Assume that $\beta_{xx}(x,t)\leq 2$
(in the sense of distributions). Then the bond call option price
$C(x,t;T_1,T_2)$ defined above is convex in $x$ at all times
$t\in[0,T_1]$.
\end{theorem}

\begin{proof}
We know from Theorem~\ref{main} (applied with $g=1$) that the price
$u$ of the $T_2$-bond is convex in $x$. 
Note that an increasing function of the bond price is
decreasing as a function of the short rate.
Thus, since the pay-off
function of the call is convex and increasing, it follows that
$C(x,T_1;T_1,T_2)=(u(x,T_1)-K)^+$ is convex. Moreover, since $u$ is
decreasing in $x$ so is $C(x,T_1;T_1,T_2)$. Since there exists a
constant $M>0$ such that \[(u(x,T_1)-K)^+\leq u(x,T_1)\leq
M\max\{e^{Mx},1\},\] 
see Corollary~\ref{cor}, the result follows
from another application of Theorem~\ref{main}.
\end{proof}

We also have a related monotonicity result. To formulate this, let
$X$ and $\tilde X$ be two diffusion processes satisfying
\[dX_t=\beta(X_t,t)\,dt+\sigma(X_t,t)\,dB_t\]
and
\[d\tilde X_t=\tilde\beta(\tilde X_t,t)\,dt+\tilde\sigma(X_t,t)\,dB_t.\]
Let $C$ and $\tilde C$ be the corresponding call option prices
written on the $T_2$-bond prices $u$ and $\tilde u$, respectively.
Then we have the following result.

\begin{theorem}
\label{callmono} Assume that $\beta(x,t)\leq \tilde \beta(x,t)$ and
$\vert\sigma(x,t)\vert\geq \vert\tilde\sigma(x,t)\vert$ for all $x$
and $t$. Also assume that either $\beta_{xx}\leq 2$ (in the 
distributional sense) at all points or
$\tilde\beta_{xx}\leq 2$ at all points. Then
\[\tilde C(x,t;T_1,T_2)\leq C(x,t;T_1,T_2)\]
for $t\leq T_1$.
\end{theorem}

\begin{proof}
First note that
\begin{equation}
\label{gtilde}
g(x):=(u(x,T_1)-K)^+\geq (\tilde u(x,T_1)-K)^+=:\tilde g(x)
\end{equation}
since $u\geq \tilde u$ by Theorem~\ref{mono}. If 
$\beta_{xx}\leq 2$ at all points, then $g$ is
decreasing and convex. It thus follows from Theorem~\ref{mono} 
and \eqref{gtilde} that
\begin{eqnarray*}
C(x,t;T_1,T_2)&=& E_{x,t}\left[\exp\left\{-\int_t^{T_1}
X_s\,ds\right\}g(X_{T_1})\right] \\ &\geq&
E_{x,t}\left[\exp\left\{-\int_t^{T_1}
\tilde X_s\,ds\right\}g(\tilde X_{T_1})\right]\\
&\geq& E_{x,t}\left[\exp\left\{-\int_t^{T_1} \tilde
X_s\,ds\right\}\tilde g(\tilde X_{T_1})\right] = \tilde
C(x,t;T_1,T_2).
\end{eqnarray*}
A similar argument can be applied if instead 
$\tilde\beta_{xx}\leq 2$ at all points. This finishes the proof.
\end{proof} 

\begin{remark}
As noted in Section~\ref{convexity}, all models in 
Table~\ref{table1} satisfy $\beta_{xx}\leq 2$. Thus 
all those models have bond call option prices which 
are convex, decreasing in the drift and increasing in the
volatility.
\end{remark}

\section{Log-convexity of bond prices}
\label{Log-convexity} 

In this section we study convexity properties
of the logarithms of bond prices. Recall that a non-negative
function $u$ is said to be log-convex if $u(\lambda
x_1+(1-\lambda)x_2)\leq u(\lambda x_1)u((1-\lambda)x_2)$ for all
$\lambda\in(0,1)$ and $x_1$, $x_2$. If $u$ is strictly positive,
then log-convexity is equivalent to the function $x\mapsto \ln u(x)$
being convex.

For simplicity, we only deal with log-convexity of bond prices, so we assume that
the pay-off function $g\equiv 1$. Thus
\begin{equation}
\label{Wg1} W(x,t)=e^{f(x)(e^{D(T-t)}-1)/D} E_{x,t}\left[\exp\left\{ -\int_t^T f(X_s)\,ds\right\}\right].
\end{equation}
Recall from Lemma~\ref{growth} that the function $W$ is
bounded on $\R\times[0,T]$. We first show that if the drift $\beta$ is
such that \eqref{driftcondition} is satisfied, then the function $W$
is also bounded away from zero.

\begin{lemma}
\label{away} Assume that $g\equiv 1$, and that
Hypothesis~\ref{hypf} and \eqref{driftcondition} hold. Then $W$
defined in \eqref{Wg1} is also bounded away from 0, i.e. there exists a constant $C>0$ such that
\[\frac{1}{C}\leq W(x,t)\leq C\]
on $\R\times[0,T]$.
\end{lemma}

\begin{proof}
Using Jensen's inequality, applied to the exponential
function and the expectation, yields
\begin{eqnarray*}
W(x,t) &=& E_{x,t}\left[\exp\left\{ \int_t^T
e^{D(s-t)}f(x)-f(X_s)\,ds \right\}\right] \\
&\geq& \exp\left\{
E_{x,t}\left[ \int_t^T e^{D(s-t)}f(x)-f(X_s)\,ds\right]\right\}\\ &=& \exp\left\{-\int_t^T e^{D(s-t)}E\left[Y_s\right]\,ds
\right\},
\end{eqnarray*}
where again $Y(s)=e^{-D(s-t)}f(X_s)-f(x)$ with initial condition
$Y_t=0$. Note that the condition \eqref{driftcondition} yields that
the drift $\tilde\beta$ of $Y$ is bounded also from above by some
constant $C$, so $E \left[Y_s\right]\leq C(s-t)$ for
$s\in[t,T]$. Consequently,
\[W(x,t) \geq \exp\left\{\int_t^T -e^{D(s-t)}C(s-t)\,ds\right\}
\geq \exp\{-Ce^{DT}T^2\}>0\] for all $t\in[0,T]$, which finishes
the proof of the lemma.
\end{proof}

\begin{theorem}
\label{log-convexity} Assume that $\alpha$ is spatially convex and
$\beta$ is spatially concave. Then the bond price $u(x,t)$ is
log-convex in the spot rate $x$.
\end{theorem}

\begin{remark}
It follows that all the models in Table~\ref{table1} have log-convex bond prices.
\end{remark}

\begin{proof}
We need to show that $\ln u(x,t)$ is spatially convex. Similar to
the proof of Theorem~\ref{main}, we approximate the function 
$\ln u$ by the function $\ln V=\ln V^f$ for some function $f$
satisfying Hypothesis~\ref{hypf}. We first assume that
Hypothesis~\ref{hyp} holds. In addition to the bounds on
the derivatives in Hypothesis~\ref{hyp}, we also assume
that 
\begin{equation}
\label{al}\vert\partial_x^k\alpha(x,t)\vert\leq C(1+x)^{1-k}\hspace{10mm} k=0,1,2
\end{equation}
for positive $x$. 

It is straightforward to check that the
function $P:=\ln V^f$ satisfies the non-linear equation
$P_\tau=\hat{\mathcal L} P$, where
\[\hat{\mathcal L} P=\alpha P_{xx}+\alpha P_x^2+\beta P_x-f,\]
with initial condition $P(x,0)=0$ (here we have again used the
parametrization in terms of time $\tau$ to maturity rather than the
physical time $t$). Moreover, there exists a constant $C$ such that
\begin{equation}
\label{vxx} \vert P_{xx}(x,\tau)\vert\leq \frac{C}{1+\vert x\vert^2}
\end{equation}
for all $(x,\tau)\in\R\times[0,T]$. Indeed, for large $\vert x\vert$
we have $f_{xx}(x)=0$, so
\[P_{xx}=(\ln V)_{xx}=(\ln W -fh)_{xx}=(\ln
W)_{xx}=\frac{WW_{xx}-W_x^2}{W^2},\] which decays like $\vert x\vert^{-2}$ according to Proposition~\ref{wdecay} and Lemma~\ref{away}.

Let $q:\R\to\R$ be a decreasing and convex function with strictly
positive second derivative such that
\[q(x)=\left\{\begin{array}{ll}
-x-(\ln \vert x\vert)^2 & \mbox{if }x\leq -K\\
-(\ln x)^2 & \mbox{if }x\geq K\end{array}\right.\] for some positive constant $K$. We now claim that there exists a
positive constant $M$ such that
\begin{equation}
\label{lett}
Mq_{xx} >\partial_x^2(\alpha q_{xx}+\beta q_x+2\alpha P_x q_x)+\vert\partial_x^2(\alpha
q_x^2)\vert
\end{equation}
at all points $(x,\tau)\in\R\times[0,T]$. 
To see this, note that $q_{xx}$ behaves like $\vert x\vert^{-2}\ln\vert x\vert$ for large $\vert x\vert$. The estimates 
\[\vert P_x+f_x(e^{D\tau}-1)/D\vert\leq \frac{C}{1+\vert
x\vert}\hspace{5mm}\mbox{ and } \hspace{5mm}\vert P_{xxx}\vert\leq
\frac{C}{1+\vert x\vert^3}\] can be obtained in the same way as
\eqref{vxx} was derived. Using these estimates,
Hypothesis~\ref{hyp} and \eqref{al}, it is straightforward to check that all terms
in
\[\partial_x^2(\alpha q_{xx}+\beta
q_x+2\alpha P_x q_x)+\vert\partial_x^2(\alpha
q_x^2)\vert\] 
decay at least like $\vert x\vert^{-2}\ln \vert x\vert$.
Consequently we can choose $M$ large so that \eqref{lett} holds.

Next, for $\ep\in(0,e^{-MT})$, define
\[P^\ep(x,\tau):=P(x,\tau)+\ep e^{M\tau}q(x).\]
Since $P_\tau=\hat{\mathcal L}P$, straightforward calculations yield that
\begin{eqnarray}
\label{persson}
\partial_x^2(P^\ep_\tau-\hat{\mathcal L} P^\ep) &=& \ep e^{M\tau} Mq_{xx}
-\ep e^{M\tau}\partial_x^2(\alpha q_{xx}+\ep e^{M\tau}\alpha
q_x^2+\beta q_x+2\alpha P_x q_x)\\
\notag
&\geq& \ep e^{M\tau}\left( Mq_{xx} -\partial_x^2(\alpha q_{xx}+\beta q_x+2\alpha P_x q_x)-\vert\partial_x^2(\alpha
q_x^2)\vert\right) >0,
\end{eqnarray}
where we used $\ep e^{M\tau}\leq 1$ and \eqref{lett}.

It follows from \eqref{vxx} and the fact that $q_{xx}$ behaves like
$x^{-2}\ln \vert x\vert$ for large $\vert x\vert$ that the set
\[\Gamma:=\{(x,\tau)\in\R\times[0,T]:P^\ep_{xx}(x,\tau)<0\}\]
is bounded. Thus, if $\Gamma$ is non-empty, then there exists a point
$(x_0,\tau_0)\in\overline \Gamma$ such that
\[\tau_0=\inf\{\tau:(x,\tau)\in\overline \Gamma \mbox{ for some }x\in\R\}.\]
Since $P^\ep_{xx}(x,0)>0$ we have $\tau_0>0$. Consequently, at the
point $(x_0,\tau_0)$ we have
\[\partial_x^2 P^\ep_\tau=\partial_\tau P^\ep_{xx} \leq 0.\]
Moreover,
\begin{eqnarray*}
\partial_x^2(\hat{\mathcal L}P^\ep) &=& \partial_x^2(
\alpha P^\ep_{xx}+\alpha (P^\ep_x)^2+\beta P^\ep_x-f)\\&=& \alpha
P^\ep_{xxxx}+(2\alpha_x+\beta)P^\ep_{xxx}+(\alpha_{xx}
+2\beta_x)P^\ep_{xx} +\alpha_{xx}(P^\ep_x)^2\\
&& +4\alpha_x P^\ep_xP^\ep_{xx}+2\alpha P^\ep_xV^\ep_{xxx}+2\alpha
(P^\ep_{xx})^2 +\beta_{xx}P^\ep_x -f_{xx}\\
&\geq& \beta_{xx}P_x^\ep\geq 0,
\end{eqnarray*}
where the first inequality is due to $P^\ep_{xxxx}\geq
P^\ep_{xxx}=P^\ep_{xx}=0$ at $(x_0,\tau_0)$, the assumption
$\alpha_{xx}\geq 0$ and the concavity of $f$, and the second
inequality follows from $P^\ep_x\leq 0$ and the assumption $\beta_{xx}\leq 0$. Thus
\begin{equation}
\label{belgium}
\partial_x^2(P^\ep_\tau-\hat{\mathcal L} P^\ep)\leq 0
\end{equation} 
at $(x_0,\tau_0)$. But this is a contradiction to \eqref{persson},
which shows that the set $\Gamma$ is empty. Thus $P^\ep$ is convex.
Letting $\ep$ tend to 0 we find that also $P$ is convex at all times
$\tau\in[0,T]$.

Now, if Hypothesis~\ref{hyp} and \eqref{al} are not satisfied, then we can
approximate $\alpha$ and $\beta$ with smooth coefficients as in the
proof of Theorem~\ref{main}. Using Proposition~\ref{approx} and then
letting $f\to x$ it is straightforward to check that also $\ln u$ is
convex, which finishes the proof.
\end{proof}

\section{Log-concavity of bond prices}
\label{Log-concavity} 
In this section we discuss concavity
properties of the logarithm $F=\ln u(x,t)$ of the bond price. Note
that $F$ satisfies the non-linear parabolic equation
\begin{equation} \label{Ueq} F_\tau=\alpha F_{xx}+\alpha F_x^2+\beta
F_x-x\end{equation} with initial condition $F(x,0)=0$. In order to
find the appropriate condition for preservation of log-concavity, we first present
a heuristic argument similar to the one presented in the beginning of
Section~\ref{convexity}.

Assume that $(x_0,\tau_0)$ is a first point where concavity is
almost lost,
i.e. $x\mapsto F(x,\tau)$ is concave for all $\tau\leq \tau_0$, and
$F_{xx}(x_0,\tau_0)=0$. Differentiating \eqref{Ueq} twice gives
\begin{eqnarray*}
\partial_\tau F_{xx} &=& \partial_x^2 F_\tau=\partial_x^2 (\alpha
F_{xx}+\alpha F_x^2+\beta F_x-x)\\
&=& \alpha F_{xxxx}+(2\alpha_x +\beta)F_{xxx}
+(\alpha_{xx}+2\beta_x) F_{xx}+\beta_{xx}F_x\\
&& + \alpha_{xx}
F_x^2+4\alpha_xF_xF_{xx}+2\alpha(F_xF_{xxx}+F_{xx}^2).
\end{eqnarray*}
Since $F_{xxxx}\leq F_{xxx}=F_{xx}=0$ at the point $(x_0,\tau_0)$ we
get
\[\partial_\tau F_{xx}\leq \beta_{xx} F_x+\alpha_{xx}F_x^2\]
at that point. Thus, since $F_x\leq 0$, we see that a sufficient
condition for preservation of log-concavity appears to be that
$\beta$ is convex and $\alpha$ concave. Since $\alpha=\sigma^2/2$ is
non-negative, however, our convention that the model is specified on
the whole real line is no longer convenient. Indeed, specifying
$\sigma$ to be 0 for negative short rates is not compatible with a
concave infinitesimal variance $\alpha$. The only possibility to
have $\alpha$ concave on the whole real line is to require it to be
constant in $x$. Such models can be shown to preserve 
log-concavity using the same methods as in Section~\ref{Log-convexity}.

\begin{theorem}
\label{log-concavity1} Assume that $\alpha$ is only time-dependent,
i.e. $\alpha(x,t)=\gamma(t)$ for some function $\gamma$. Also assume
that $\beta$ is convex in $x$ for each fixed time $t$. Then the
$T$-bond price $u(x,t)$ is log-concave in $x$ at every fixed time
$t\in[0,T]$.
\end{theorem}

\begin{proof}
The proof follows along the same lines as
Theorem~\ref{log-convexity} with some minor modifications. The
function $f$ of Hypothesis~\ref{hypf} needs to be replaced by a
convex function which equals $x$ for positive $x$ and is
constant for $x\leq -K$, where $K$ is some positive constant.
With this new $f$, Lemma~\ref{away} and Proposition~\ref{wdecay} remain valid if we modify $\beta$ to be
linear for $x$ large positive. 
\end{proof}

The conditions in Theorem~\ref{log-concavity1} are only satisfied by
the models V and HW in Table~\ref{table1}. To investigate the
remaining models, we are forced to leave the tractable setting of
diffusions defined on the whole real line and instead
consider models specified on a half-line. Such models, however,
typically lead to partial differential equations with degenerate
coefficients. 

For simplicity, we assume that $X$ is defined on the positive real 
axis $[0,\infty)$. We will also assume that $\alpha(0,t)=0$ and 
$\beta(0,t)\geq 0$. Note that under these conditions, no 
boundary behavior of the process at $x=0$ needs to be specified.
Let 
\begin{equation}
\label{w} 
w(x,t)=e^{h(t)x}u(x,t),
\end{equation}
where $h(t)=(e^{D(T-t)}-1)/D$. 
By arguing as in the proofs of Lemma~\ref{growth} and
Lemma~\ref{away}, it can be shown that if 
$\beta-Dx$ and $\alpha$ are bounded for large $x$,
then there exists a positive constant $C$ such that
\begin{equation}
\label{waway}
C^{-1}\leq w(x,t)\leq C
\end{equation}
for all $(x,t)\in[0,\infty)\times [0,T]$.
To apply the techniques used in previous sections, we also need 
estimates of the derivatives of the function $w$. 
We have the following result.

\begin{lemma}
\label{derivatives}
Assume that $X$ is specified on the positive real axis with $\alpha(0,t)=0$, $\beta(0,t)\geq 0$, $\alpha(x,t)>0$ for $x\in
(0,\infty)$, that 
$\alpha$ and $\beta$ are smooth and that $\alpha$ and 
$\beta-Dx$ are constant in $x$ for large $x$.
Then there exists a constant $C>0$ such that for $k=0,1,2,...$ we have
\begin{equation}
\label{craig}
\vert\partial_x^k w(x,t)\vert\leq \frac{C}{1+x^k}
\end{equation}
for $(x,t)\in(0,\infty)\times[0,T]$.
\end{lemma}

\begin{proof} 
As noted above, the case $k=0$ can be proven analogously to
Lemma~\ref{growth}. It is straightforward to check
that $\hat w:=w-1$ is the bounded classical solution to 
\begin{equation} 
\label{what}
\left\{\begin{array}{ll}
\hat w_t+\alpha \hat w_{xx}+\hat\beta \hat w_x+\gamma 
\hat w+\gamma=0 & t<T\\
w=0 & t=T,\end{array}\right.
\end{equation}
where 
\[\hat\beta=\beta-2\alpha h\hspace{7mm}\mbox{ and }\hspace{7mm} \gamma=\alpha h^2+(Dx-\beta)h.\]
Note that
\begin{equation}
\label{bo}
\vert\partial_x^k \hat\beta(x,t)\vert\leq C(1+x)^{1-k}
\hspace{7mm}\mbox{ and }\hspace{7mm} \vert\partial_x^k\gamma(x,t)\vert \leq C(1+x)^{-k}
\end{equation}
for some constant $C$. By stochastic representation,
\[\hat w(x,t)=E\left[\int_t^T\exp\left\{
\int_t^s\gamma(Y^{x,t}_r,r)\,dr\right\}\gamma(Y^{x,t}_s)\, ds\right],\]
where $Y^{x,t}$ is a diffusion given by
\[dY^{y,t}_s=\hat\beta(Y^{x,t}_s,s)+\sigma(Y^{x,t}_s,s)\,dB_s\]
and the indices indicate that $X^{x,t}_t=x$.
As in the proof of Lemma~\ref{wdecay}, the bounds in \eqref{bo} and
Gronwall's lemma can be applied to prove that $\hat w_x$ is
bounded and satisfies $\hat w_x(x,T)=0$. Thus $\hat w_x$ is a
bounded solution to the equation obtained by differenting 
the equation \eqref{what}, and using the maximum principle 
the estimate \eqref{craig} can be established for $k=1$.
The rest of the proof follows inductively by 
treating the differentiated equation as above.
\end{proof}

\begin{remark}
It follows from Lemma~\ref{derivatives} that the derivatives
$\partial_x^k\partial_t^l u$, $k+2l\leq 4$, are continuous
up to the spatial boundary $x=0$.
\end{remark}

\begin{theorem}
\label{log-concavity2} Assume that $X$ is specified on the 
positive real axis with $\alpha(0,t)=0$ and $\beta(0,t)\geq 0$.
Also assume that $\alpha$ is concave in $x$ and $\beta$ is
convex in $x$ at any fixed time $t\in[0,T]$.
Then the bond price $u(x,t)$ is log-concave in $x$.
\end{theorem}

\begin{proof}
We will assume that $\alpha$ and $\beta$
satisfy the conditions in Lemma~\ref{derivatives}. We also
assume that there exists a constant $\eta>0$ such that
$\alpha(x,t)=\gamma(t)x$ and $\beta_{xx}(x,t)=0$ for $x\leq \eta$ and for some function
$\gamma(t)\geq \eta$. The general case follows by approximation.

Define the function $F(x,\tau)=\ln u(x,\tau)$, where $u$ is as in
\eqref{u}. As above, $F$ satisfies the non-linear equation
\[F_\tau=\hat{\mathcal L} F\]
where
\[\hat{\mathcal L}F=\alpha F_{xx}+\alpha F_x^2+\beta F_x-x.\]
It follows from \eqref{waway} and Lemma~\ref{derivatives} that
\[\vert F_x(x,t)+h(t)\vert\leq C(1+x)^{-1},\]
\[\vert F_{xx}(x,t)\vert\leq C(1+x)^{-2}\]
and
\[\vert F_{xxx}(x,t)\vert\leq C(1+x)^{-3}.\]
Let $q:(0,\infty)\to\R$ be a smooth, increasing and concave function
with strictly negative second derivative such that
\[q(x)=\left\{\begin{array}{ll}x-x^2 &\mbox{if }x<1/C\\
(\ln x)^2 &\mbox{if } x>C\end{array}\right.\] for some constant
$C>0$. We claim that there is a constant $M>0$ so large that
\begin{equation}
\label{litau}
Mq_{xx} <\partial_x^2(\alpha q_{xx}+\beta q_x+2\alpha F_x q_x)-\vert\partial_x^2(\alpha q_x^2)\vert
\end{equation}
at all points $(x,\tau)\in[0,\infty)\times[0,T]$.
Indeed, the right hand side is bounded for small $x$, and
for large values of $x$ the right hand side decays at least as
fast as $x^{-2}\ln x$. Consequently, $M$ can
be chosen so that \eqref{litau} holds at all points. 

Now, for $\ep\in(0,e^{-MT})$, define
\[F^\ep(x,\tau)=F(x,\tau)+\ep e^{M\tau}q(x).\]
Then
\begin{eqnarray}
\label{M} 
\partial_x^2(F_\tau^\ep-\hat{\mathcal L} F^\ep) &=&
\ep e^{M\tau} \left( Mq_{xx} -\partial_x^2(\alpha q_{xx}+\ep
e^{M\tau}\alpha q_x^2+\beta q_x+2\alpha F_x q_x)\right)\\
\notag &\leq& \ep e^{M\tau} \left( Mq_{xx} -\partial_x^2(\alpha q_{xx}+\beta q_x+2\alpha F_x q_x)+\vert\partial_x^2 (\alpha q_x^2)\vert\right) <0
\end{eqnarray}
according to \eqref{litau}. Next, define the set
\[\Gamma=\{(x,t)\in\R\times [0,T]: F^\ep_{xx}(x,t)<0\}.\]
Since $F_{xx}$ decays at least like $x^{-2}$ we have that $\Gamma\subseteq [0,R)\times [0,T]$
for some $R$, so $\overline \Gamma$ is compact. Let 
$(x_0,\tau_0)$ be a point such that
\[\tau_0=\inf\{\tau>0: (x,\tau)\in\Gamma
\mbox{ for some }x\in[0,\infty)\}.\]
If $x_0>0$, then arguing as before the inequality
\eqref{belgium}, it is straightforward to check that 
\[\partial_x^2(F_\tau^\ep-\hat{\mathcal L} F^\ep)\geq 0\]
at $(x_0,\tau_0)$, which contradicts \eqref{M}. Therefore,
assume that $x_0=0$. Then, at the point $(x_0,\tau_0)$ we have
\begin{equation}
\label{nuder}
\partial_x^2 F_\tau^\ep=\partial_\tau F_{xx}^\ep\geq 0
\end{equation}
and
\begin{eqnarray}
\label{sahlin}
\partial_x^2(\hat{\mathcal L} F^\ep) &=& \partial_x^2
(\alpha F_{xx}^\ep+\alpha (F_x^\ep)^2+\beta F_x^\ep-x)\\
\notag &=& 
\alpha_{xx}F_{xx}^\ep+2\alpha_xF_{xxx}^\ep+\alpha F_{xxxx}^\ep
+\alpha_{xx}F_{xx}^\ep+4\alpha_x F_x^\ep F_{xx}^\ep \\
\notag && 
+2\alpha (F_x^\ep F_{xxx}^\ep+(F_{xx}^\ep)^2)+\beta_{xx}F_x^\ep
+2\beta_x F_{xx}^\ep+\beta F_{xxx}^\ep\leq 0,
\end{eqnarray}
where we have used $\alpha=0$, $\beta\geq 0$, $F_{xx}=0$ and $F_{xxx}\leq 0$ at the point $(x_0,\tau_0)$. The inequalities 
\eqref{M}, \eqref{sahlin} and \eqref{nuder} form a contradiction.
This shows that the set $\Gamma$ is empty, so $F^\ep$ is
convex at all times. Letting $\ep\to 0$, it follows that also $F$ is
convex in $x$.
\end{proof}

\begin{remark}
From the results in Sections~\ref{Log-convexity} and
\ref{Log-concavity} we find that if $\alpha$ and $\beta$ are both
concave and convex, i.e. linear, then both log-convexity and
log-concavity are preserved. In our terminology such models thus
give rise to log-affine bond prices. Of course, these models are
usually referred to as having an {\it affine term structure} and
they play an important rôle in interest rate theory. Our sufficient
conditions for the existence of an affine term structure are
well-known, see for instance Chapter~17 of \cite{B}. However, the
results of Sections~\ref{Log-convexity} and \ref{Log-concavity}
offer a background to these seemingly ad hoc
conditions.
\end{remark}

\section{Conclusions}
\label{summary}
In this paper we have conducted a study of convexity of solutions to
the term structure equation. We show that if the drift $\beta$
satisfies $\beta_{xx}\leq 2$, then the bond prices are convex in the
current short rate, increasing in the volatility of the short rate
and decreasing in the drift. Similar results hold for call options
written on bond prices. For models with regular coefficients, the
condition $\beta_{xx}\leq 2$ is also a necessary condition for
preservation of convexity. We also have a general comparison
theorem: if a model has smaller drift and larger volatility than
another model, and at least one of them has a drift satisfying the
condition above, then the first model has the larger bond prices.
For bond call options the analogous result holds.

We also study convexity properties of the logarithm of a bond price
corresponding to the relative sensitivity of the bond price to
changes in the short rate. This relative sensitivity is often
described by the duration, i.e. the negative of the derivative of
the logarithm. We show that if the drift $\beta$ is
concave and the square $\sigma^2$ of the volatility is convex, then
bond prices are log-convex (a decreasing duration). Similarly, if $\beta$
is convex and $\sigma^2$ is concave, then bond prices are
log-concave (an increasing duration). We also note that if we demand that the
price is log-convex {\it and} log-concave, we recover the well-known
sufficient conditions for a model to admit an affine term structure.

Our findings for some commonly used models are summarized in Table~\ref{table2} below.

\begin{table}[ht!]
\begin{center}
\begin{tabular}{|c|c|c|c|c|c|c|}
\hline Model & Dynamics  & AB & AO & C & LCV & LCC\\
\hline \hline V & $dX=k(\theta-X)\,dt+\sigma\,dB$  & Yes & Yes & Yes & Yes & Yes\\
\hline CIR & $dX=k(\theta-X)\,dt+\sigma\sqrt X\,dB$  & Yes & Yes & Yes & Yes & Yes\\
\hline D & $dX=bX\,dt+\sigma X\,dB$  &Yes & No & Yes & Yes & No\\
\hline EV & $dX=X(\eta-a\ln X)\,dt+\sigma X\,dB$  & No & No & Yes & Yes & No\\
\hline HW & $dX=k(\theta_t-X)\,dt+\sigma\,dB$  & Yes & Yes & Yes & Yes & Yes\\
\hline BK & $dX=X(\eta_t-a\ln X)\,dt+\sigma X\,dB$  & No & No & Yes & Yes & No\\
\hline MM & $dX=X(\eta_t-(\lambda-\frac{\gamma}{1+\gamma t})\ln
X)\,dt+\sigma X\,dB$  & No & No & Yes & Yes & No\\
\hline
\end{tabular}\caption{\label{table2}
In the three last columns it is indicated
which models preserve convexity for option prices and bond call options, 
log-convexity of bond prices and log-concavity of bond prices,
respectively.}
\end{center}
\end{table}


\begin{thebibliography}{99999}

\bibitem{A1}
Alvarez, L.H.R. On the properties of $r$-excessive mappings for a
class of diffusions. Ann. Appl. Probab. 13 (2003) 1517-1533.

\bibitem{A}
Alvarez, L.H.R. On the form and risk-sensitivity of zero coupon
bonds for a class of interest rate models. Insurance Math. Econom.
28 (2001) 83-90.

\bibitem{BMO}
Bahlali, K., Mezerdi, B. and Ouknine, Y. Pathwise uniqueness and
approximation of solutions of stochastic differential equations.
Séminaire de Probabilités, XXXII, 166-187, Lecture Notes in Math.,
1686, Springer, Berlin, 1998.

\bibitem{BJ}
Bellamy, N. and Jeanblanc, M. Incompleteness of markets driven by a mixed diffusion. Finance Stoch. 4 (2000), 209-222.

\bibitem{BR}
Bergenthum, J., and Ruschendorf, L. Comparison of option prices in
semimartingale models. Finance Stoch. 10 (2006), 222-249.

\bibitem{BGW}
Bergman, Y., Grundy, B. and Wiener, Z. General properties of option
prices. J. Finance 51 (1996) 1573-1610.

\bibitem{B}
Bj\"ork, T. Arbitrage Theory in Continuous Time, Oxford University
Press, New York, (1998).

\bibitem{BS}
Bj\"ork, T. and Svensson, L. On the existence of finite-dimensional realizations for nonlinear forward rate models. Math. Finance 11  (2001), 205-243.

\bibitem{B1}
Borell, C. Geometric inequalitites in option pricing, in: Convex
Geometric Analysis, in: Math. Sci. Res. Inst. Publ., vol. 34,
Cambridge Univ. Press, Cambridge, 1999, 29-51.

\bibitem{B2}
Borell, C. Isoperimetry, log-concavity, and elasticity of option
prices, in: P. Wilmott, H. Rasmussen (Eds.), New Directions in
Mathematical Finance, Wiley, 2002, 73-91.

\bibitem{BM}
Brigo, D. and Mercurio, F. Interest rate models - theory and
practice. Springer Finance. Springer-Verlag, Berlin, 2001.

\bibitem{E}
Ekstr\"om, E. Properties of American option prices. Stochastic
Process. Appl. 114 (2004), 265-278.

\bibitem{EJT}
Ekstr\"om, E., Janson, S. and Tysk, J. Superreplication of options
on several underlying assets. J. Appl. Probab. 42 (2005), 27-38.

\bibitem{ET1}
Ekstr\"om, E. and Tysk, J. Properties of option prices in models
with jumps. (2005) To appear in Math. Finance.

\bibitem{ET2}
Ekstr\"om, E. and Tysk, J. Convexity preserving jump-diffusion
models for option pricing. To appear in J. Math. Anal. Appl. (2006).

\bibitem{ET3}
Ekstr\"om, E. and Tysk, J. The American put is log-concave in the
log-price. J. Math. Anal. Appl. 314 (2006), 710-723.

\bibitem{EKJPS}
El Karoui, N., Jeanblanc-Picqué, M. and Shreve, S.  Robustness of
the Black and Scholes formula. Math. Finance 8 (1998), no. 2,
93-126.

\bibitem{ER}
Eriksson, J. Monotonicity in the volatility of single-barrier
options, Int. J. Theor. Appl. Finance 9 (2006) 987-996.

\bibitem{H}
Hobson, D. Volatility misspecification, option pricing and
superreplication via coupling. Ann. Appl. Probab. 8 (1998)
193-205.

\bibitem{JT1}
Janson, S. and Tysk, J. Volatility time and properties of option
prices. Ann. Appl. Probab. 13 (2003) 890-913.

\bibitem{JT2}
Janson, S. and Tysk, J. Preservation of convexity of solutions to
parabolic equations. J. Differential Equations 206 (2004)
182-226.

\bibitem{K}
Kolesnikov, A.V. On diffusion semigroups preserving the
log-concavity, J. Funct. Anal. 186 (2001) 196-205.

\bibitem{L}
Lieberman, G. M. Initial regularity for solutions of degenerate
equations. Nonlinear Anal. 14 (1990), 525-536.

\bibitem{M}
Merton, R. Theory of rational option pricing. Bell J. Econom. and
Management Sci. 4 (1973), 141-183.

\bibitem{RY}
Revuz, D. and Yor, M. Continuous martingales and Brownian motion. Third edition. Grundlehren der Mathematischen Wissenschaften, 293. Springer-Verlag, Berlin, 1999.

\bibitem{Y}
Yong, J. Remarks on some short rate term structure models. J. Ind.
Manag. Optim. 2 (2006) 119-134.

\end{thebibliography}
\end{document}